
%
%
%

\input amstex
\loadbold
\documentstyle{amsppt}
\NoBlackBoxes

\pagewidth{32pc}
\pageheight{44pc}
\magnification=\magstep1

\def\ZZ{\Bbb{Z}}

\def\RR{\Bbb{R}}
\def\CC{\Bbb{C}}
\def\PP{\Bbb{P}}
\def\OO{\Cal{O}}
\def\Spec{\operatorname{Spec}}
\def\Sym{\operatorname{Sym}}
\def\Sing{\operatorname{Sing}}
\def\zeros{\operatorname{div}}
\def\rank{\operatorname{rk}}
\def\tr{\operatorname{tr}}
\def\Ker{\operatorname{Ker}}
\def\JHGr{\operatorname{Gr}^{JH}}
\def\twedge{\bar{\bigwedge}}
\def\rest#1#2{\left.{#1}\right\vert_{{#2}}}
\def\achern#1#2{\widehat{c}_{#1}(#2)}

\topmatter
\title
Arithmetic Bogomolov-Gieseker's inequality
\endtitle
\rightheadtext{}
\author Atsushi Moriwaki \endauthor
\leftheadtext{}
\address
Department of Mathematics, Faculty of Science,
Kyoto University, Kyoto, 606-01, Japan
\endaddress
\curraddr
Department of Mathematics, University of California,
Los Angeles, 405 Hilgard Avenue, Los Angeles, California 90024, USA
\endcurraddr
\email moriwaki\@math.ucla.edu \endemail
\date July, 1993 \enddate
\thanks
THIS IS A TENTATIVE VERSION.
\endthanks
\abstract
Let $f : X \to \Spec(\ZZ)$ be an arithmetic variety of dimension $d \geq 2$
and $(H, k)$ an arithmetically ample Hermitian line bundle on $X$, that is,
a Hermitian line bundle with the following properties:
\roster
\item "(1)" $H$ is $f$-ample.

\item "(2)" The Chern form $c_1(H_{\infty}, k)$ gives a K\"{a}hler form
            on $X_{\infty}$.

\item "(3)" For every irreducible horizontal subvariety $Y$
            (i.e. $Y$ is flat over $\Spec(\ZZ)$), the height
            $\achern{1}{\rest{(H, k)}{Y}}^{\dim Y}$ of $Y$ is positive.
\endroster
Let $(E, h)$ be a rank $r$ vector bundle on $X$.
In this paper, we will prove that
if $E_{\infty}$ is semistable with respect to $H_{\infty}$ on
each connected component of $X_{\infty}$, then
$$
\left\{ \achern{2}{E, h} - \frac{r-1}{2r} \achern{1}{E, h}^2 \right\} \cdot
\achern{1}{H, k}^{d-2} \geq 0.
$$
Moreover, if the equality of the above inequality holds,
then $E_{\infty}$ is projectively flat and
$h$ is a weakly Einstein-Hermitian metric.
\endabstract
\toc\nofrills{Table of Contents}
\widestnumber\head{10}
\specialhead {} Introduction \endspecialhead
\head 1. Intersection pair and positive cone \endhead
\head 2. Unstability theorem \endhead
\head 3. Restriction theorem \endhead
\head 4. Sufficiently ample divisor \endhead
\head 5. Strictly effective section \endhead
\head 6. Proof of Main Theorem \endhead
\specialhead {} References \endspecialhead
\endtoc
\endtopmatter
\vfil\eject

\document

\subhead Introduction
\endsubhead
Let $f : X \to \Spec(\ZZ)$ be a $d$-dimensional arithmetic
variety (i.e. $X$ is a $d$-dimensional integral scheme and
$f$ is a projective and flat morphism with the smooth generic fiber).
Let $(E, h)$ be a rank $r$ Hermitian vector bundle on $X$.
In \cite{Mo}, we proved that, if $d=2$ and
$E_{\infty}$ $(= E \otimes_{\ZZ} \CC)$ is semistable
on $X_{\infty}$ $(= X \otimes_{\ZZ} \CC)$, then
$$
  \achern{2}{E, h} - \frac{r-1}{2r} \achern{1}{E, h}^2 \geq 0.
$$
We would like to generalize the above inequality to a higher
dimensional arithmetic variety. For this purpose, let
$(H, k)$ be a Hermitian line bundle on $X$ such that
$H$ is $f$-ample and the Chern form $c_1(H_{\infty}, k)$
gives a K\"{a}hler form on $X_{\infty}$.
Unfortunately, even if $E_{\infty}$
is $H_{\infty}$-semistable,
$$
\left\{ \achern{2}{E, h} - \frac{r-1}{2r} \achern{1}{E, h}^2 \right\} \cdot
\achern{1}{H, k}^{d-2}
$$
is not non-negative in general. For example, if
$E_{\infty}$ is not projectively flat and
$c$ is a sufficiently large positive number, then
$$
\left\{ \achern{2}{E, h} - \frac{r-1}{2r} \achern{1}{E, h}^2 \right\} \cdot
\achern{1}{H, ck}^{d-2} < 0.
$$
This indicates us that we need a good condition for $(H, k)$. That is
``arithmetical ampleness'' due to S. Zhang \cite{Zh},
which is, roughly speaking, a natural arithmetic analogy of Nakai-Moishezon
criterion for ampleness. More precisely,
$(H, k)$ is said to be {\it arithmetically ample} if the following
conditions are satisfied.
\roster
\item "(1)" $H$ is $f$-ample.

\item "(2)" The Chern form $c_1(H_{\infty}, k)$ gives a K\"{a}hler form
            on $X_{\infty}$.

\item "(3)" For every irreducible horizontal subvariety $Y$
            (i.e. $Y$ is flat over $\Spec(\ZZ)$), the height
            $\achern{1}{\rest{(H, k)}{Y}}^{\dim Y}$ of $Y$ is positive.
\endroster
Using an arithmetically ample Hermitian line bundle,
our main theorem of this paper is the following.

\proclaim{Main Theorem}
Let $f : X \to \Spec(\ZZ)$ be an arithmetic variety of dimension $d \geq 2$
and $(H, k)$ an arithmetically ample Hermitian line bundle on $X$.
Let $(E, h)$ be a rank $r$ vector bundle on $X$.
If $E_{\infty}$ is semistable with respect to $H_{\infty}$ on
each connected component of $X_{\infty}$, then we have
$$
\left\{ \achern{2}{E, h} - \frac{r-1}{2r} \achern{1}{E, h}^2 \right\} \cdot
\achern{1}{H, k}^{d-2} \geq 0.
$$
Moreover, if the equality of the above inequality holds,
then $E_{\infty}$ is projectively flat and
$h$ is a weakly Einstein-Hermitian metric.
\endproclaim

Here, we would like to explain technical aspects of the proof of Main Theorem.
In the geometric case, the Bogomolov-Gieseker's inequality
is derived from the Bogomolov-Gieseker's inequality on surfaces and
Mehta-Ramanathan's restriction theorem \cite{MR}.
But the arithmetic case is not so simple.
For example, let $\phi$ be a section of $H^0(X, H)$. Then
$$
\multline
\left\{ \achern{2}{E, h} - \frac{r-1}{2r} \achern{1}{E, h}^2 \right\} \cdot
\achern{1}{H, k}^{d-2} = \\
      \left\{ \achern{2}{\rest{(E, h)}{\zeros(\phi)}} -
              \frac{r-1}{2r} \achern{1}{\rest{(E, h)}{\zeros(\phi)} }^2
              \right\} \cdot
       \achern{1}{\rest{(H, k)}{\zeros(\phi)}}^{d-3} \\
 - \int_{X_{\infty}} \log(\sqrt{k(\phi, \phi)})
         \left\{ c_2(E_{\infty}, h) - \frac{r-1}{2r} c_1(E_{\infty}, h)^2
\right\}
         c_1(H_{\infty}, k)^{d-3}.
\endmultline
$$
To proceed induction of $\dim X$, we need a section $\phi$ with
the following properties.
\roster
\item "(a)" $\zeros(\phi_{\infty})$ is smooth.

\item "(b)" $\rest{E_{\infty}}{\zeros(\phi_{\infty})}$ is semistable.

\item "(c)" $||\phi||_{\sup} < 1$.
\endroster
By S. Zhang's result \cite{Zh} concerning the existence of strictly effective
sections, we can find a section $\phi$ satisfying
(a) and (c), replacing $H$ by $H^m$ if necessarily (cf. Corollary~1.3 and
Theorem~2.3).
Unfortunately, Mehta-Ramanathan's restriction theorem does not guarantee
a section with (a), (b) and (c). We need a more powerful restriction theorem,
that is,

\proclaim{Bogomolov's restriction theorem}{\rm (cf. Theorem~3.1 and
Lemma~4.3.2)}
Let $X$ be a smooth projective variety of dimension $d \geq 2$
over an algebraically closed field of characteristic zero and
$H$ an ample divisor on $X$.
Let $E$ be a semistable torsion free sheaf on $X$ with respect to $H$.
Then, there are an effectively determined integer $m_0$ and
closed points $x_1, \ldots, x_s$ of $X$ such that, for all $m \geq m_0$,
if $C \in |mH|$ is normal and $C \cap \{ x_1, \ldots, x_s \} = \emptyset$,
then $\rest{E}{C}$ is semistable with respect to $\rest{H}{C}$.
\endproclaim

\noindent
The original result of Bogomolov \cite{Bo2} was restricted to the surface case,
but in this paper, we generalize it to a higher dimensional variety.

\medskip
In \S1, \S2 and \S3, we will prove a generalization of
Bogomolov's restriction theorem.
\S4 and \S5 are preliminaries for the proof of the main theorem.
\S6 is devoted to the proof of the main theorem.

\subhead 1. Intersection pair and positive cone \endsubhead
Let $X$ be a smooth projective variety of dimension $d \geq 2$ and
$H_1, \ldots, H_{d-2}$ ample divisors on $X$.
We set $\Cal{H} = (H_1, \ldots, H_{d-2})$ and
$N^1(X) = \left( \{ \hbox{divisors on $X$} \}/\equiv \right) \otimes \RR$.
Here we define a natural pairing
$$
  (\ \cdot \ )_\Cal{H} : N^1(X) \times N^1(X) \to \RR
$$
by $(x \cdot y)_\Cal{H} = (x \cdot y \cdot H_1 \cdots H_{d-2})$.
The set $\{ x \in N^1(X) \mid (x \cdot x)_\Cal{H} > 0 \}$
consists of two connected components. One of them contains all ample
divisors. This component is called the positive cone of $X$ and
is denoted by $P(X; \Cal{H})$.
The following lemma is useful for later purpose.

\proclaim{Lemma 1.1}
Let $\overline{P}(X; \Cal{H})$ be the topological closure
of $P(X; \Cal{H})$. Then, we have the following.
\roster
\item "(1)" $(x \cdot y)_\Cal{H} >  0$ for all $x \in P(X; \Cal{H})$
            and $y \in \overline{P}(X; \Cal{H}) \setminus \{ 0 \}$.

\item "(2)" $x \in P(X; \Cal{H})$ if and only if
            $(x \cdot y)_\Cal{H} > 0$ for all
            $y \in \overline{P}(X; \Cal{H}) \setminus \{ 0 \}$.
\endroster
\endproclaim

\demo{Proof}
(1) By Hodge index theorem, we have a basis
$\{ e_1, \ldots, e_n \}$ of $N^1(X)$ such that
$e_1 \in P(X; \Cal{H})$, $(e_1 \cdot e_1)_\Cal{H} = 1$,
$(e_2 \cdot e_2)_\Cal{H} = (e_3 \cdot e_3)_\Cal{H} = \cdots =
(e_n \cdot e_n)_\Cal{H} = -1$ and that
$(e_i \cdot e_j)_\Cal{H} = 0$ for $i \not= j$.
We set $x = x_1 e_1 + \cdots + x_n e_n$ and
$y = y_1 e_1 + \cdots + y_n e_n$.
Then,
$$
x_1 > \sqrt{(x_2)^2 + \cdots + (x_n)^2}\quad\hbox{and}\quad
y_1 \geq \sqrt{(y_2)^2 + \cdots + (y_n)^2}.
$$
Since $y \not= 0$, $y_1 > 0$.
Thus, by Schwartz' inequality, we get
$$
\align
 (x \cdot y)_\Cal{H} & = x_1 y_1 - x_2 y_2 - \cdots - x_n y_n \\
                     & \geq x_1 y_1 - \sqrt{(x_2)^2 + \cdots + (x_n)^2}
                                      \sqrt{(y_2)^2 + \cdots + (y_n)^2} \\
                     & = y_1 (x_1 - \sqrt{(x_2)^2 + \cdots + (x_n)^2}) + \\
                     & \quad \sqrt{(x_2)^2 + \cdots + (x_n)^2}
                             (y_1 - \sqrt{(y_2)^2 + \cdots + (y_n)^2}) \\
                     & > 0.
\endalign
$$

(2) By (1), it is sufficient to show that
$x \in P(X; \Cal{H})$ if $(x \cdot y)_\Cal{H} > 0$ for all
$y \in \overline{P}(X; \Cal{H}) \setminus \{ 0 \}$.
We set $x = x_1 e_1 + \cdots + x_n e_n$.
Since $x_1 = (x \cdot e_1)_\Cal{H} > 0$, we may assume that $x_1 = 1$.
Moreover, we may assume that $(x_2)^2 + \cdots + (x_n)^2 \not= 0$.
We set
$$
       y = e_1 + \frac{1}{\sqrt{(x_2)^2 + \cdots + (x_n)^2}}
                 (x_2 e_2 + \cdots + x_n e_n).
$$
Then $y \in \overline{P}(X; \Cal{H}) \setminus \{ 0 \}$.
Thus,
$$
   (x \cdot y)_\Cal{H} = 1 - \sqrt{(x_2)^2 + \cdots + (x_n)^2} >  0.
$$
Therefore, we have
$$
  (x \cdot x)_\Cal{H} = 1 - \left( (x_2)^2 + \cdots + (x_n)^2 \right) > 0.
\eqno{\qed}
$$
\enddemo

\subhead 2. Unstability theorem \endsubhead
First of all, we will introduce several notations.
Let $X$ be a smooth projective variety of dimension $d \geq 2$ and
$H, H_1, \ldots, H_{d-2}$ ample divisors on $X$.
We set $\Cal{H} = (H_1, \ldots, H_{d-2})$.
For a torsion free sheaf $E$ on $X$, an averaged degree
$\mu(E; H, H_1, \ldots, H_{d-2})$ of $E$
with respect to $H, H_1, \ldots, H_{d-2}$ and a discriminant
$\delta_\Cal{H}(E)$ of $E$ with respect to $\Cal{H}$
are defined by
$$
 \mu(E; H, H_1, \ldots, H_{d-2}) =
 \frac{(c_1(E) \cdot H \cdot H_1 \cdots H_{d-2})}{\rank E},
$$
$$
 \delta_\Cal{H}(E) =
 \left( \left(
 \frac{\rank E - 1}{2 \rank E} c_1(E)^2 - c_2(E)
 \right) \cdot H_1 \cdots H_{d-2} \right).
$$
We say $E$ is stable (resp. semistable) with respect to
$H, H_1, \ldots, H_{d-2}$ if,
for all subsheaves $F$ with $0 \subsetneq F \subsetneq E$,
$$
    \mu(F; H, H_1, \ldots, H_{d-2}) <
    \mu(E; H, H_1, \ldots, H_{d-2}).
$$
$$
    \hbox{(resp. $\mu(F; H, H_1, \ldots, H_{d-2}) \leq
                  \mu(E; H, H_1, \ldots, H_{d-2})$)}
$$
Moreover, for torsion free sheaves $E$ and $F$ on $X$,
we set
$$
  d(F, E) = \frac{c_1(F)}{\rank F} - \frac{c_1(E)}{\rank E}.
$$

Let $0 \to S \to E \to Q \to 0$ be an exact sequence of
torsion free sheaves on $X$. We can easily see that
$$
 \delta_\Cal{H}(E) = \delta_\Cal{H}(S) + \delta_\Cal{H}(Q) +
             \frac{(\rank E)(\rank S)}{2 \rank Q} (d(S, E)^2)_\Cal{H}.
\tag 2.1
$$
In particular, if $(d(S, E)^2)_\Cal{H} \geq 0$, then
$$
 \delta_\Cal{H}(E) \leq \delta_\Cal{H}(S) + \delta_\Cal{H}(Q) +
                \frac{\rank E(\rank E - 1)}{2} (d(S, E)^2)_\Cal{H}.
\tag 2.2
$$

The purpose of this section is to give a generalization of
Bogomolov's unstability theorem \cite{Bo1} to a higher dimensional
projective variety.

\proclaim{Theorem 2.3}
Let $X$ be a smooth projective variety of dimension $d \geq 2$
over an algebraically closed field of characteristic zero
and $\Cal{H} = (H_1, \ldots, H_{d-2})$ a sequence of
ample divisors on $X$.
Let $E$ be a torsion free sheaf on $X$. If $\delta_\Cal{H}(E) > 0$,
there is a saturated subsheaf $F$ of $E$ with
$d(F, E) \in P(X; \Cal{H})$.
\endproclaim

\demo{Proof}
Let $H$ be another ample divisor on $X$.
We set
$$
   W(G) =
   \{ x \in \overline{P}(X; \Cal{H}) \setminus \{ 0 \}
      \mid (d(G, E) \cdot x)_\Cal{H} > 0 \}
$$
for a saturated subsheaf $G$ of $E$ with $(d(G, E) \cdot H)_\Cal{H} > 0$.
Then, by virtue of (2) of Lemma~1.1, $d(G, E) \in P(X; \Cal{H})$
if and only if $W(G) = \overline{P}(X; \Cal{H}) \setminus \{ 0 \}$.
Here we claim:

\proclaim{Claim 2.3.1} The set
$$
\{ [d(G, F)] \in N^1(X) \mid \hbox{$G$ is a saturated subsheaf of $E$ with
                         $(d(G, E) \cdot H)_\Cal{H} > 0$} \}
$$
is finite.
\endproclaim

For this purpose, it is sufficient to show the following lemma.

\proclaim{Lemma~2.3.2}
Let $T$ be a torsion free sheaf on $X$ and $M$ a real number. Then,
the set
$$
\{ [c_1(L)] \in N^1(X) \mid \hbox{$L$ is a rank $1$ subsheaf of $T$ with
$(c_1(L) \cdot H)_\Cal{H} \geq M$} \}
$$
is finite.
\endproclaim

\demo{Proof}
It is easy to see that there is a filtration of $T$:
$0 = T_0 \subset T_1 \subset \cdots \subset T_{l-1} \subset T_l = T$
such that $T_{i}/T_{i-1}$ is a rank $1$ torsion free sheaf for every $i$.
Let $L_i$ be double dual of $T_{i}/T_{i-1}$.
Let $L$ be a rank $1$ subsheaf of $T$ with
$(c_1(L) \cdot H)_\Cal{H} \geq M$.
Pick up $i$ with $L \not\subset T_{i-1}$ and $L \subset T_i$.
Then, since $L \to T_{i}/T_{i-1}$ is non-trivial,
there is an effective divisor $D_L$ on $X$ such that
$c_1(L) +  D_L = c_1(L_i)$. Thus
$$
 (D_L \cdot H)_\Cal{H} \leq
 (c_1(L_i) \cdot H)_\Cal{H} - M.
$$
Therefore $D_L$ has a bounded degree.
It follows that $D_L$ sits in a bounded family of effective divisors on $X$.
Hence, we have our lemma.
\qed
\enddemo

Since $\delta_\Cal{H}(E) > 0$, by Corollary~4.7 of \cite{Mi},
$E$ is not semistable with respect to $H, H_1, \ldots, H_{d-2}$.
Thus, there is a saturated subsheaf $F$ of $E$ with
$(d(F, E) \cdot H)_\Cal{H} > 0$.
Then, by (2.1),
$$
 \delta_\Cal{H}(E) = \delta_\Cal{H}(F) + \delta_\Cal{H}(E/F) +
             \frac{(\rank E)(\rank F)}{2 \rank (E/F)}
             (d(F, E)^2)_\Cal{H}.
\tag 2.3.3
$$

First, we consider the case where $\rank E = 2$.
Since $\rank F = \rank E/F = 1$, we have
$\delta_\Cal{H}(F) \leq 0$ and $\delta_\Cal{H}(E/F) \leq 0$.
It follows $(d(F, E)^2)_\Cal{H} > 0$ by (2.3.3).

\medskip
In general, we prove this theorem by induction on $\rank E$.
Here we claim that

\proclaim{Claim 2.3.4}
If $(d(F, E)^2)_\Cal{H} \leq 0$, then there is a saturated subsheaf
$F_1$ of $E$ such that $(d(F_1, E) \cdot H)_\Cal{H} > 0$ and
$W(F) \subsetneq W(F_1)$.
\endproclaim

Since $(d(F, E)^2)_\Cal{H} \leq 0$, by (2.3.3), we have either
$\delta_\Cal{H}(F) > 0$ or $\delta_\Cal{H}(E/F) > 0$.

If $\delta_\Cal{H}(F) > 0$, then by hypothesis of induction
there is a saturated subsheaf $F_1$ of $F$ with
$d(F_1, F) \in P(X; \Cal{H})$. Here
since $d(F_1, E) = d(F_1, F) + d(F, E)$,
we have $W(F) \subsetneq W(F_1)$.

If $\delta_\Cal{H}(E/F) > 0$, then by hypothesis of induction
there is a saturated subsheaf $F_1$ of $E$ such that
$F \subset F_1$ and $d(F_1/F, E/F) \in P(X; \Cal{H})$. Here by
an easy calculation, we get
$$
d(F_1, E) = \frac{\rank (F_1/F) }{\rank F_1} d(F_1/F, E/F) +
             \frac{\rank F \rank (E/F_1)}
                  {\rank F_1 \rank (E/F) } d(F, E).
$$
Therefore, $W(F) \subsetneq W(F_1)$. Thus we have our claim

\medskip
We set $F_0 = F$. If $(d(F_0, E)^2)_\Cal{H} > 0$, $F_0$ is our
desired subsheaf. Otherwise, by Claim~2.3.4, there is a saturated
subsheaf $F_1$ of $E$ such that $(d(F_1, E) \cdot H)_\Cal{H} > 0$
and $W(F_0) \subsetneq W(F_1)$. If $(d(F_1, E)^2)_\Cal{H} > 0$,
then we have our theorem. Otherwise, by Claim~2.3.4, we get
a saturated subsheaf $F_2$ of $E$ such that
$(d(F_2, E) \cdot H)_\Cal{H} > 0$ and $W(F_1) \subsetneq W(F_2)$.
Here we assume that continuing these procedures, we can not get
a saturated subsheaf $F_n$ with $d(F_n, E) \in P(X; \Cal{H})$.
Then, there is a sequence
$\{ F_0, F_1, F_2, \ldots, F_n, \ldots \}$ of
saturated subsheaves of $E$ such that $(d(F_i, E) \cdot H))_\Cal{H} >  0$
for all $i \geq 0$ and
$$
W(F_0) \subsetneq W(F_1) \subsetneq \cdots \subsetneq
W(F_n) \subsetneq \cdots \subsetneq \overline{P}(X; \Cal{H}) \setminus \{ 0 \}.
$$
In particular, the numerical classes $[d(F_i, E)]$ of $d(F_i, E)$ are
distinct. This contradicts to Claim~2.3.1.
\qed
\enddemo

\proclaim{Corollary 2.4}
Let $X, H_1, \ldots, H_{d-2}$ be same as in Theorem~{\rm 2.3}.
Let $E$ be a torsion free sheaf on $X$.
If $\delta_\Cal{H}(E) > 0$, then we have the following.
\roster
 \item "(1)" There is a saturated subsheaf $F$ of $E$ with
             $\delta_\Cal{H}(F) \leq 0$ and
             $d(F, E) \in P(X; \Cal{H})$.

 \item "(2)" There is a saturated subsheaf $T$ of $E$ with
             $\delta_\Cal{H}(E/T)\leq 0$ and
             $d(T, E) \in P(X; \Cal{H})$.
\endroster
\endproclaim

\demo{Proof}
(1) We set
$$
  D(E) = \{ F \mid
       \hbox{$F$ is a saturated subsheaf of $E$ with
             $d(F, E) \in P(X; \Cal{H})$} \}.
$$
By Theorem~2.3, $D(E) \not= \emptyset$.
Moreover, by virtue of Claim~2.3.1,
the image of $D(E)$ in $N^1(X)$ is finite.
Therefore, there is an element $F$ of $D(E)$ such that
$(d(F, E)^2)_\Cal{H}$ is maximal.
Let us see $\delta_\Cal{H}(F) \leq 0$.
If $\delta_\Cal{H}(F) > 0$,
by Theorem~2.3, there is a saturated subsheaf $L$ of $F$
with $d(L, F) \in P(X; \Cal{H})$.
Then, $d(L, E) = d(L, F) + d(F, E) \in P(X; \Cal{H})$.
Thus, by (1) of Lemma~1.1,
$$
 (d(L, E)^2)_\Cal{H} = \left((d(L, F) + d(F, E))^2 \right)_\Cal{H} >
                            (d(F, E)^2)_\Cal{H}.
$$
This is a contradiction.

(2) Applying (1) to
the dual $E^{\vee}$ of $E$, we have a saturated subsheaf $T'$ of
$E^{\vee\vee}$ with $\delta_\Cal{H}(E^{\vee\vee}/T') \leq 0$ and
$d(T', E^{\vee\vee})_\Cal{H} \in P_{X, \Cal{H}}$.
We set $T = T' \cap E$. Then we have
$d(T, E) = d(T', E^{\vee\vee})$ and
$\delta_\Cal{H}(E/T) \leq \delta_\Cal{H}(E^{\vee\vee}/T') \leq 0$.
Thus we have the second assertion.
\qed
\enddemo

\proclaim{Corollary 2.5}
Let $X, H_1, \ldots, H_{d-2}$ be same as in Theorem~{\rm 2.3}.
Let $E$ be a torsion free sheaf on $X$.
If $\delta(E)_\Cal{H} > 0$, there is a saturated subsheaf $L$ of $E$
such that $\delta_\Cal{H}(L) \leq 0$,
$d(L, E) \in P(X; \Cal{H})$ and that
$$
     \frac{\rank E (\rank E - 1)}{2} (d(L, E)^2)_\Cal{H} \geq
     \delta(E)_\Cal{H}.
$$
\endproclaim

\demo{Proof}
By (2) of Corollary~2.4, there is a filtration of $E$:
$ 0 = T_0 \subset T_1 \subset \cdots \subset T_{l-1} \subset T_l = E$
of length $l \geq 2$ with the following properties:
  \roster
   \item "(a)" $T_i/T_{i-1}$ is torsion free for every $1 \leq i \leq l$.

   \item "(b)" $\delta_\Cal{H}(T_i/T_{i-1}) \leq 0$
               for every $1 \leq i \leq l$

   \item "(c)" $d(T_{i-1}, T_{i}) \in P(X; \Cal{H})$
               for every $2 \leq i \leq l$
  \endroster
Thus by (1) of Lemma~1.1 and (2.2), we get
$$
\align
\delta_\Cal{H}(E)
& \leq \sum_{i=1}^{l} \delta_\Cal{H}(T_{i}/T_{i-1}) +
       \sum_{i=2}^{l} \frac{\rank T_{i} (\rank T_{i} - 1)}{2}
                      (d(T_{i-1}, T_{i})^2)_\Cal{H} \\
          & \leq \frac{\rank E (\rank E - 1)}{2} \sum_{i=2}^{l}
                      (d(T_{i-1}, T_{i})^2)_\Cal{H} \\
          & \leq \frac{\rank E (\rank E - 1)}{2}
                 \left( \left( \sum_{i=2}^{l} d(T_{i-1}, T_{i})
                 \right)^2 \right)_\Cal{H} \\
          & = \frac{\rank E (\rank E - 1)}{2} (d(T_{1}, E)^2)_\Cal{H}
\endalign
$$
Therefore, $T_1$ is our desired subsheaf.
\qed
\enddemo

\subhead 3. Restriction theorem
\endsubhead
Let $E$ be a rank $r$ vector bundle on a smooth projective variety of
dimension $d$ over an algebraically closed field of characteristic zero.
We assume that $E$ is semistable with respect to an ample divisor $H$.
V. Mehta and A. Ramanathan \cite{MR} proved that,
for a sufficiently large $m$ and
a general member $C$ of $|mH|$, the restriction $\rest{E}{C}$ to $C$
is also semistable. Unfortunately, by their method, we can not
find an effective estimation of $m$. For example, this is very important
to see boundedness of a family of semistable vector bundles.
In \cite{Fl}, H. Flenner found an effective estimation of $m$.
More precisely, he proved that if $m$ satisfies an inequality
$$
    \frac{\displaystyle \binom{d+m}{m} - m - 1 }{m} >
    (H^d) \max \left\{ \frac{r^2 - 1}{4}, 1 \right\},
$$
then, for a general hypersurface $C$ in $|mH|$, the restriction
$\rest{E}{C}$ to $C$ is also semistable.

Next it is very natural to ask whether
semistability is preserved by a special restriction or not.
Recently, F. A. Bogomolov \cite{Bo2} gives an answer of the above question
for the case where $d = 2$. His result says
that the restriction of a stable vector bundle to a smooth member of $|mH|$
is also stable. In the semistable case, his result doesn't not hold
in general. For example, there is a rank $2$ vector bundle $E$ on $\PP^2$
with an exact sequence $0 \to \OO_{\PP^2} \to E \to m_x \to 0$,
where $m_x$ is the maximal ideal at a point $x \in \PP^2$.
Then, $E$ is semistable, but $\rest{E}{C}$ is not semistable
for all curves $C$ passing through $x$.
This example shows us that we must take care of
pinching points of the Jordan-H\"{o}lder filtration in the semistable case.
In this section, we would like to give a generalization
of Bogomolov's restriction theorem \cite{Bo2} to a higher dimensional variety
including a semistable case.
To state the main theorem of this section, first of all,
we will introduce a Jordan-H\"{o}lder filtration of a semistable vector bundle.

If a vector bundle $E$ is semistable with respect to $H, H_1, \ldots, H_{d-2}$,
there is a filtration of $E$:
$$
  0 = E_0 \subset E_1 \subset \cdots \subset E_{l-1} \subset E_l = E
$$
with the following properties:
\roster
\item "(1)" $E_i/E_{i-1}$ is torsion free for every $1 \leq i \leq l$.

\item "(2)" $E_i/E_{i-1}$ is stable with respect to
$H, H_1, \ldots, H_{d-2}$ for every $1 \leq i \leq l$.

\item "(3)" $\mu(E_i/E_{i-1}; H, H_1, \ldots, H_{d-2}) =
             \mu(E; H, H_1, \ldots, H_{d-2})$ for every $1 \leq i \leq l$.
\endroster
The above filtration is called a Jordan-H\"{o}lder filtration of $E$.
It is well know that $\bigoplus_{i=1}^{l} E_i/E_{i-1}$ does not
depend on the choice of a filtration. So we denote this sheaf
by $\JHGr(E)$.

The main theorem of this section is the following:

\proclaim{Theorem 3.1}
Let $X$ be a smooth projective variety of dimension $d \geq 2$
over an algebraically closed field of characteristic zero and
$H, H_1, \ldots, H_{d-2}$ ample divisors on $X$.
We set $\Cal{H} = (H_1, \ldots, H_{d-2})$.
Let $E$ be a semistable torsion free sheaf with respect to
$H, H_1, \ldots, H_{d-2}$ on $X$.
Let $m$ be a positive integer with
$$
 m > \max_{0 < p < \rank E} \{ -2 \rank(\twedge^p E)
                                  \delta_\Cal{H}(\twedge^p E) \}
$$
and $C$ a divisor in $|mH|$ such that $C$ is normal and
$\rest{\JHGr(\twedge^p E)}{C}$ has no torsion for all $0 < p < \rank E$,
where $\twedge^p$ means $p$-th exterior power modulo torsion.
Then, $\rest{E}{C}$ is semistable with respect to
$\rest{H_1}{C}, \ldots, \rest{H_{d-2}}{C}$.
Moreover, if $E$ is reflexive, then
$$
 \max_{0 < p < \rank E} \{ -2 \rank(\twedge^p E)
                              \delta_\Cal{H}(\twedge^p E) \} =
  -2 \binom{\rank E}{[\rank E/2]} \binom{\rank E - 2}{[\rank E/2] - 1}
     \delta_\Cal{H}(E).
$$
{\rm (}Note that $\twedge^p E$ is semistable with respect to
$H, H_1, \ldots, H_{d-2}$ for every $0 < p < \rank E$.{\rm )}
\endproclaim

The following lemma is a key for the proof of Theorem~3.1.

\proclaim{Lemma 3.2}
Let $X$ be a smooth projective variety of dimension $d \geq 2$
over an algebraically closed field of characteristic zero
and $H, H_1, \cdots, H_{d-2}$ ample divisors on $X$.
We set $\Cal{H} = (H_1, \ldots, H_{d-2})$.
Let $E$ be a rank $r$ stable torsion free sheaf
with respect to $H, H_1, \ldots, H_{d-2}$.
Let $m$ be a positive integer with $m > -2 r\delta_\Cal{H}(E)$ and
$C$ a divisor in $|mH|$ such that $C$ is normal and
$\rest{E}{C}$ has no torsion.
Then, for all rank $1$ torsion free quotient sheaves $Q$ of $\rest{E}{C}$,
we have
$$
\mu(\rest{E}{C}; \rest{H_1}{C}, \ldots, \rest{H_{d-2}}{C}) <
\mu(Q; \rest{H_1}{C}, \ldots, \rest{H_{d-2}}{C}).
$$
\endproclaim

\demo{Proof}
Assume that there is a rank $1$ torsion free quotient sheaf
$Q$ of $\rest{E}{C}$ with
$$
   \mu(\rest{E}{C}; \rest{H_1}{C}, \ldots, \rest{H_{d-2}}{C}) \geq
   \mu(Q; \rest{H_1}{C}, \ldots, \rest{H_{d-2}}{C}).
$$
We set $F = \Ker(E \to Q)$. Then it is easy to see that
$c_1(F) = c_1(E) - C$ and
$$
\align
   (c_2(F) \cdot H_1 \cdots H_{d-2}) & =
   (c_2(E) \cdot H_1 \cdots H_{d-2}) -
   \deg(\rest{E}{C}; \rest{H_1}{C}, \ldots, \rest{H_{d-2}}{C}) \\
   & \qquad + \deg(Q; \rest{H_1}{C}, \ldots, \rest{H_{d-2}}{C}).
\endalign
$$
Thus we have
$$
\align
\delta_\Cal{H}(F) & = \delta_\Cal{H}(E) + \frac{r-1}{2r}(C^2)_\Cal{H} \\
& \qquad + \mu(\rest{E}{C}; \rest{H_1}{C}, \ldots, \rest{H_{d-2}}{C}) -
           \mu(Q; \rest{H_1}{C}, \ldots, \rest{H_{d-2}}{C}) \\
& \geq \delta_\Cal{H}(E) + \frac{m^2(r-1)}{2r} (H^2)_\Cal{H}.
\endalign
$$
Since $-2r\delta_\Cal{H}(E)$ is a non-negative integer,
it is easy to see that
$$
      m > \sqrt{\frac{-2r\delta_\Cal{H}(E)}{(r-1)(H^2)}}.
$$
It follows $\delta_\Cal{H}(F) > 0$.
Thus, by Corollary~2.5, there is a saturated subsheaf $L$ of $F$ such that
$d(L, F) \in P(X; \Cal{H})$ and
${\displaystyle \frac{r(r-1)}{2} (d(L, F)^2)_\Cal{H} \geq
\delta_\Cal{H}(F)}$.
In particular, we get
$$
\frac{r(r-1)}{2} \left( \left( d(L, E) + \frac{mH}{r} \right)^2
\right)_\Cal{H} \geq
\delta_\Cal{H}(E) + \frac{m^2(r-1)}{2r}(H^2)_\Cal{H} ,
$$
which implies that
$$
\frac{r(r-1)}{2} \left( d(L, E) \cdot \left( d(L, E) + \frac{2m H}{r} \right)
\right)_\Cal{H} \geq \delta_\Cal{H}(E).
$$

Here we claim that
$$
    \left( d(L, E) \cdot \left( d(L, E) + \frac{m H}{r} \right)
    \right)_\Cal{H} < 0.
$$
If $(d(L, E)^2)_\Cal{H} \leq 0$,
the assertion is trivial because $E$ is stable.
So we may assume that $(d(L, E)^2)_\Cal{H} > 0$.
Then, $-d(L, E) \in P(X; \Cal{H})$.
On the other hand, $d(L, F) = d(L, E) + (m/r)H \in P(X; \Cal{H})$.
Thus we have our assertion by (1) of Lemma~1.1.

By the above claim, we get
$$
\delta_\Cal{H}(E) < \frac{m(r-1)}{2} ( d(L, E) \cdot H)_\Cal{H}.
$$
Let $l$ be a rank of $L$.
Since $lr(d(L, E) \cdot H)_\Cal{H}$ is a negative integer,
we have $(d(L, E) \cdot H)_\Cal{H} \leq -1/lr$. Therefore
$$
 \delta_\Cal{H}(E) < \frac{-m(r-1)}{2lr} \leq \frac{-m}{2r}.
$$
Thus $m < -2r\delta_\Cal{H}(E)$. This is a contradiction.
\qed
\enddemo

\proclaim{Corollary 3.3}
Let $X$ be a smooth projective variety of dimension $d \geq 2$
over an algebraically closed field of characteristic zero
and $H, H_1, \cdots, H_{d-2}$ ample divisors on $X$.
We set $\Cal{H} = (H_1, \ldots, H_{d-2})$.
Let $E$ be a semistable torsion free sheaf
with respect to $H, H_1, \ldots, H_{d-2}$.
Let $m$ be a positive integer with $m > -2\rank(E)\delta_\Cal{H}(E)$ and
$C$ a divisor in $|mH|$ such that $C$ is normal and
$\rest{\JHGr(E)}{C}$ has no torsion.
Then, for all rank $1$ torsion free quotient sheaves $Q$ of $\rest{E}{C}$,
we have
$$
\mu(\rest{E}{C}; \rest{H_1}{C}, \ldots, \rest{H_{d-2}}{C}) \leq
\mu(Q; \rest{H_1}{C}, \ldots, \rest{H_{d-2}}{C}).
$$
\endproclaim

\demo{Proof}
Let $0 = E_0 \subset E_1 \subset \cdots \subset E_{l-1} \subset E_l = E$
be a Jordan-H\"{o}lder filtration of $E$.
We set $Q_i = E_i/E_{i-1}$.
By (2.1) and Hodge index theorem,
$\delta_\Cal{H}(E) \leq \sum_{i=1}^l \delta_\Cal{H}(Q_i)$.
Therefore, we have $-2\rank(E)\delta_\Cal{H}(E) \geq
-2\rank(Q_i)\delta_\Cal{H}(Q_i)$ for all $i$.
Let $j$ be the minimal number such that $\rest{E_j}{C} \to Q$
is non-trivial. Then, we have a non-trivial homomorphism
$\rest{Q_j}{C} \to Q$. Therefore, by Lemma~3.2, we get
$$
\mu(\rest{Q_j}{C}; \rest{H_1}{C}, \ldots, \rest{H_{d-2}}{C}) \leq
\mu(Q; \rest{H_1}{C}, \ldots, \rest{H_{d-2}}{C}).
$$
Thus we have our Corollary.
\qed
\enddemo

\demo{Proof of Theorem~{\rm 3.1}}
Let us start the proof of Theorem~3.1.
Let $Q$ be a rank $p$ torsion free quotient sheaf of $\rest{E}{C}$.
Then, $\twedge^p Q$ is a rank $1$ torsion free quotient sheaf of
$\twedge^p \left(\rest{E}{C}\right)$.
Since $\bigwedge^p \left(\rest{E}{C}\right) \simeq
\rest{\left(\bigwedge^p E\right)}{C}$,
we have $\twedge^p \left(\rest{E}{C}\right) \simeq
\rest{\left(\bigwedge^p E\right)}{C}/\hbox{torsion}$.
On the other hand,
$\rest{\left(\bigwedge^p E \right)}{C}/\hbox{torsion} \simeq
\rest{\left(\twedge^p E\right)}{C}$
because there is a surjective homomorphism
$\rest{\left(\bigwedge^p E\right)}{C} \to
\rest{\left(\twedge^p E\right)}{C}$ and
$\rest{\left(\twedge^p E\right)}{C}$ is torsion free.
Therefore, $\twedge^p Q$ is a rank 1 torsion free quotient sheaf
of $\rest{\left(\twedge^p E\right)}{C}$. Thus, by Corollary~3.3, we get
$$
\mu(\rest{\left(\twedge^p E\right)}{C}; \rest{H_1}{C}, \ldots,
\rest{H_{d-2}}{C}) \leq
\mu(\twedge^p Q; \rest{H_1}{C}, \ldots, \rest{H_{d-2}}{C}),
$$
which implies that
$$
\mu(\rest{E}{C}; \rest{H_1}{C}, \ldots, \rest{H_{d-2}}{C}) \leq
\mu(Q; \rest{H_1}{C}, \ldots, \rest{H_{d-2}}{C}).
$$
Therefore, $\rest{E}{C}$ is semistable.

If $E$ is reflexive, then, by a calculation of Chern classes, we have
$$
\delta_\Cal{H}(\twedge^p E) = \binom{\rank E - 2}{p - 1}\delta_\Cal{H}(E).
$$
(For this calculation, we may assume that $d = 2$ and $E$ is locally free.)
Therefore,
$$
 \max_{0 < p < \rank E} \{ -2 \rank(\twedge^p E)
                              \delta_\Cal{H}(\twedge^p E) \} =
  -2 \binom{\rank E}{[\rank E/2]} \binom{\rank E - 2}{[\rank E/2] - 1}
     \delta_\Cal{H}(E).
$$
Thus, we get the last assertion of Theorem~3.1
\qed
\enddemo

\subhead 4. Sufficiently ample divisor
\endsubhead
In this section, we will consider
an estimation of the degree of the locus of singular divisors
in a complete linear system. Let $X$ be a smooth projective scheme
over an algebraically closed field.
A divisor $H$ is said to be {\it sufficiently ample} if,
for all $x \not= y \in X$, $\OO_X(H) \otimes m_x^2 m_y$ is generated by
global sections, where $m_x$ and $m_y$ is the maximal ideals at $x$ and $y$.
If $H$ is sufficiently ample, then it is easy to see that,
for all $x \in X$, $\OO_X(H) \otimes m_x$ is generated by
global sections. Thus, $H$ is very ample.
Conversely, let try to see that a higher multiple of an ample divisor
is sufficiently ample.

\proclaim{Lemma 4.1}
Let $X$ be a smooth projective scheme over an algebraically closed field and
$H$ an ample divisor on $X$.
Then, there is a positive integer $m_0$ such that,
if $m \geq m_0$, then $mH$ is sufficiently ample.
\endproclaim

\demo{Proof}
We consider two closed subschemes $\Delta_1$ and $\Delta_2$ in
$X \times X \times X$ given by
$$
    \Delta_1 = \{ (x, y, z) \mid x = y \} \quad\hbox{and}\quad
    \Delta_2 = \{ (x, y, z) \mid x = z \}.
$$
Let $I_{\Delta_1}$ and $I_{\Delta_2}$ be the defining ideals of
$\Delta_1$ and $\Delta_2$.
Let $p_{23} : X \times X \times X \to X \times X$ and
$p_1 : X \times X \times X \to X$ be the natural projections to
the second-third factor and the first factor.
Since $p_1^*(\OO_X(H))$ is $p_{23}$-ample, there is a positive integer
$m_0$ such that, if $m \geq m_0$, then
$R^i{p_{23}}_*(p_1^*(\OO_X(mH)) \otimes I_{\Delta_1} I_{\Delta_2}) = 0$
for all $i > 0$ and
$$
  {p_{23}}^*{p_{23}}_*(p_1^*(\OO_X(mH)) \otimes I_{\Delta_1} I_{\Delta_2}))
  \to p_1^*(\OO_X(mH)) \otimes I_{\Delta_1} I_{\Delta_2}
$$
is surjective.
Thus, we have our lemma.
\qed
\enddemo

Here we introduce several notations. Let $X$ be a smooth
projective scheme over an algebraically closed field and
$X = X_1 \cup \cdots \cup X_l$ a decomposition into connected
components. Let $H$ be a divisor on $X$. We set
$$
  |H| = \PP(H^0(X_1, \OO_{X_1}(H))) \times \cdots \times
  \PP(H^0(X_l, \OO_{X_l}(H))).
$$
(Note that, if $X$ is not connected, then $|H|$ does not coincide
with the usual complete linear system $\PP(H^0(X, \OO_X(H)))$.)
A hypersurface $Z$ in $|H|$ is said to be {\it decomposable} if
there are hypersurfaces $Z_i$ in $\PP(H^0(X_i, \OO_{X_i}(H)))$
such that $Z = p_1^*(Z_1) + \cdots + p_l^*(Z_l)$, where
$p_i$'s are the natural projections $|H| \to \PP(H^0(X_i, \OO_{X_i}(H)))$.
It is easy to see that
if a hypersurface $Z$ is decomposable, then there are homogeneous polynomials
$f_i \in \Sym^{k_i}((H^0(X_i, \OO_{X_i}(H)))^{\vee})$ such that
$Z = \zeros(f_1 \cdots f_l)$.
Moreover, we denote by $\Sing(|H|)$ the set of
all singular divisors in $|H|$.
Using these notations, we have the following theorem.

\proclaim{Theorem 4.2}
Let $X$ be a smooth projective scheme of equi-dimension $d$
{\rm (}i.e. every connected component is of dimension $d${\rm )} over
an algebraically closed field and
$H$ a sufficiently ample divisor on $X$.
Then, $\Sing(|H|)$ is a decomposable hypersurface in $|H|$ of degree
$$
    \sum_{i=0}^d (i+1)(c_{d-i}(\Omega_X^1) \cdot H^i).
$$
Moreover, if $X$ is connected, then $\Sing(|H|)$ is irreducible.
\endproclaim

\demo{Proof}
First of all, we need the following lemma.

\proclaim{Lemma 4.2.1}
Let $X$ and $H$ be as in Theorem~{\rm 4.2}. If $X$ is connected, then
$\Sing(|H|)$ is an irreducible hypersurface in $|H|$ and
there is a non-empty Zariski open set $U$ of $\Sing(|H|)$ such that
if $D \in U$, then $D$ has only one isolated ordinary double point.
\endproclaim

\demo{Proof}
We set
$$
\Sigma = \{ (x, D) \in X \times |H| \mid x \in D \}.
$$
Let $p : \Sigma \to X$ and $q : \Sigma \to |H|$ be the natural projections.
Moreover, we set
$$
B = \{ (x, D) \in X \times |H| \mid
\hbox{$x \in D$ and $D$ is singular at $x$.} \}.
$$
We denote by $B_x$ the fiber of $p : B \to X$ at $x$.
Since $H$ is very ample, as in the proof of \cite{Ha, Theorem~II.8.18},
$B_x = \PP(H^0(X, \OO_X(H) \otimes m_x^2))$ and
$\dim B_x = n - d - 1$, where $n = \dim |H|$.
Thus, $B$ is irreducible and of dimension $n - 1$.
Therefore, since $\Sing(|H|) = p(B)$, $\Sing(|H|)$ is irreducible and
of dimension $\leq n - 1$.
Hence, in order to see our lemma,
it is sufficient to see that, for all $x \in X$,
there is a non-empty Zariski open set $U_x$ of $B_x$ such that
if $D \in U_x$, then $D \setminus \{ x \}$ is smooth and
$D$ has an ordinary double point at $x$.
Moreover, since
$$
     H^0(X, \OO_X(H) \otimes m_x^2) \to \OO(H) \otimes m_x^2/m_x^3
$$
is surjective, it is sufficient to find
an open set $U'_x$ of $B_x$ such that if $D \in U'_x$, then
$D \setminus \{ x \}$ is smooth.
For this purpose, we consider the following scheme.
$$
T = \{ (y, D) \in (X \setminus \{ x \}) \times |H| \mid
       \hbox{$x, y \in D$ and $D$ is singular at $x$ and $y$.} \}.
$$
Clearly we have $T \cap p^{-1}(y) = \PP(H^0(X, \OO_X(H) \otimes m_x^2 m_y^2))$.
Since
$$
   H^0(X, \OO_X(H) \otimes m_x^2) \to \OO_X(H) \otimes (\OO_X/m_y)
$$
and
$$
   H^0(X, \OO_X(H) \otimes m_x^2 m_y) \to \OO_X(H) \otimes (m_y/m_y^2)
$$
are surjective, we get
$\dim (T \cap p^{-1}(y)) = n - 2d - 2$.
Therefore, $\dim T = n - d - 2$.
Thus $T$ is a proper closed subset of $B_x$.
Therefore, we have our assertion.
\qed
\enddemo

Let us go back to the proof of Theorem~4.2.
First, we consider the case where $X$ is connected.
Then, by Lemma~4.2.1,
there is a pencil $\{ H_{\lambda} \}_{\lambda \in \PP^1}$ on $X$
with the following properties.
\roster
\item "(1)" $\bigcap_{\lambda \in \PP^1} H_{\lambda}$ is smooth.

\item "(2)" If $l$ is a corresponding line of the pencil
            $\bigcap_{\lambda \in \PP^1} H_{\lambda}$ in
            $|H|$, then
            $$\#(\Sing(|H|) \cap l) = \deg(\Sing(|H|)).$$

\item "(3)" $H_{\lambda}$ has at most one ordinary double point
            for all $\lambda \in \PP^1$.
\endroster
Let $\mu : Y \to X$ be a blowing-up along
$\bigcap_{\lambda \in \PP^1} H_{\lambda}$
and $f : Y \to \PP^1$ the induced morphism.
Let
$$
   0 \to f^*(\Omega_{\PP^1}^1) \to \Omega_Y^1 \to \Omega_{X/\PP^1}^1 \to 0.
$$
the canonical exact sequence.
Then, by (2), (3) and \cite{Fu, Example~3.2.16},
$$
 \deg(\Sing(|H|)) = c_d(\Omega_Y^1 \otimes (f^*(\Omega_{\PP^1}^1))^{\vee}).
$$
Therefore, by \cite{Fu, Example~3.2.2 and Theorem~15.4}, we have
$$
 \align
  \deg(\Sing(|H|)) & = c_d(\Omega_Y^1) + 2 c_{d-1}(\Omega_{H_{\eta}}^1) \\
    & = c_d(\Omega_X^1) + c_{d-2}(\Omega_C^1) + 2 c_{d-1}(\Omega_{H_{\eta}}^1)
 \endalign
$$
where $H_{\eta}$ is a general element of
$\{ H_{\lambda} \}_{\lambda \in \PP^1}$ and
$C = \bigcap_{\lambda \in \PP^1} H_{\lambda}$.
On the other hand,
$$
c_{d-1}(\Omega_{H_{\eta}}^1) = \sum_{i=1}^d (c_{d-i}(\Omega_X^1) \cdot H^i)
\quad\hbox{and}\quad
c_{d-2}(\Omega_{C}^1) = \sum_{i=2}^d (i-1)(c_{d-i}(\Omega_X^1) \cdot H^i).
$$
Hence, we have
$$
\deg(\Sing(|H|)) = \sum_{i=0}^d (i+1)(c_{d-i}(\Omega_X^1) \cdot H^i).
$$

Next we consider a general case.
Let $X = X_1 \cup \cdots \cup X_l$ be a decomposition into
connected components.
Then,
$$
\Sing(|H|) = \bigcup_{j=1}^{l} |\rest{H}{X_1}| \times \cdots \times
             \Sing(|\rest{H}{X_j}|) \times \cdots \times |\rest{H}{X_l}|.
$$
Therefore, $\Sing(|H|)$ is a decomposable hypersurface and
$$
\deg(\Sing(|H|)) = \sum_{j=1}^l \deg(\Sing(|\rest{H}{X_j}|)).
$$
Thus,
$$
\deg(\Sing(|H|)) =
\sum_{j=1}^l \sum_{i=0}^d (i+1)(c_{d-i}(\Omega_{X_j}^1) \cdot H^i)
= \sum_{i=0}^d (i+1)(c_{d-i}(\Omega_X^1) \cdot H^i).
$$
Therefore, we have our theorem.
\qed
\enddemo

\proclaim{Corollary 4.3}
Let $X$ be a smooth projective scheme of equi-dimension $d \geq 2$
over an algebraically closed field of characteristic zero and
$H$ ample divisors on $X$.
Let $E$ be a torsion free sheaf on $X$ such that
$E$ is semistable with respect to $H$ on each connected component.
Then, there are closed points $x_1, \ldots, x_s$ of $X$ such that,
if $m$ is sufficiently large and $C$ is a smooth divisor of $|mH|$ with
$C \cap \{ x_1 , \ldots, x_s \} = \emptyset$, then
$\rest{E}{C}$ is semistable with respect to $\rest{H}{C}$.
Moreover, let $Z_m$ be the set of all divisors $D$ of $|mH|$ such that
$D$ is singular or $D \cap \{ x_1 , \ldots, x_s \} \not= \emptyset$.
Then, if $m$ is sufficiently large, $Z_m$ is a decomposable
hypersurface of $|mH|$ at most of degree
$$
    \sum_{i=0}^d (i+1)(c_{d-i}(\Omega_X^1) \cdot H^i)m^i + s.
$$
\endproclaim

\proclaim{Lemma 4.3.1}
Let $A$ be a commutative ring with the identity,
$M$ a $A$-module and $N$ a $A$-submodule of $M$.
If $a$ is not a zero-divisor for $M/N$, then $N/aN \to M/aM$ is injective.
\endproclaim

\demo{Proof}
Assume that $N/aN \to M/aM$ is not injective. Then, there is an element
$x \in N$ such that $x \not\in aN$ and $x \in aM$.
Thus, we have $y \in M$ with $x = a y$.
Since $y \not\in N$, $y \not\equiv 0 \mod N$. On the other hand,
$ay \equiv 0 \mod N$. Thus, $a$ is a zero-divisor for $M/N$.
This is a contradiction.
\qed
\enddemo

\proclaim{Lemma 4.3.2}
Let $X$ be a smooth scheme over an algebraically closed field
and $F$ a torsion free sheaf on $X$.
Then, there are closed points $\{ x_1, \ldots, x_s \}$ of $X$
such that, for any smooth divisor $D$ on $X$,
if $D \cap \{ x_1 , \ldots, x_s \} = \emptyset$,
$\rest{F}{D}$ has no torsion.
\endproclaim

\demo{Proof}
Let $E$ be the double dual of $F$.
Since $E$ is locally a second syzygy sheaf, by Lemma~4.3.1,
$\rest{E}{D}$ has no torsion for all smooth divisors $D$.
Let $\{ P_1, \ldots, P_s \}$ be the set of associated primes of $E/F$ and
$V_i = \Spec(\OO_X/P_i)$.
Pick up closed points $x_i \in V_i$ with $x_i \not\in V_j$ for $j \not= i$.
Let $D$ be a smooth divisor with $\{ x_1, \cdots, x_s \} \cap D = \emptyset$.
Then, $D \not\in P_i$ for all $i$.
Thus, by Lemma~3.1.1, $\rest{F}{D} \to \rest{E}{D}$ is injective.
Therefore, $\rest{F}{D}$ has no torsion.
\qed
\enddemo

\demo{Proof of Corollary~{\rm 4.3}}
Let us start of the proof of Corollary~4.3.
The first assertion is an immediate consequence of
Theorem~3.1 and Lemma~4.3.2.
Let $x$ be a close point of $X$. Then, the set of all divisors in $|mH|$
passing through $x$ is a decomposable hypersurface in $|mH|$ of degree $1$.
Hence the second assertion is obtained by Lemma~4.1 and Theorem~4.2.
\qed
\enddemo

\subhead 5. Strictly effective section
\endsubhead
In this section, we will consider strictly effective sections
of a Hermitian line bundle on an arithmetic variety. First of all,
we introduce a simple notation.
If an integer $a$ has the prime factorization
$a = \pm p_1^{e_1} \cdot p_2^{e_2} \cdots p_r^{e_r}$, then
we set
$\operatorname{rad}(a) = p_1 \cdot p_2 \cdots p_r$.

\proclaim{Lemma 5.1}
Let $f : X \to \Spec(\ZZ)$ be an arithmetic variety and $L$ an $f$-ample
line bundle. Let $x_1, x_2, \ldots, x_r$ be distinct points of $X$
such that the residue field $\kappa(x_i)$ at $x_i$ has the positive
characteristic
for every $i$. We set
$l = \operatorname{rad}(
\operatorname{char}(\kappa(x_1))\cdots\operatorname{char}(\kappa(x_r)))$.
Then, there is a positive integer $n_0$ such that, for all $n \geq n_0$,
if $e_1, e_2, \ldots, e_{p(n)}$ are generators of $H^0(X, L^n)$ as
a $\ZZ$-module, there are integers $a_1, \ldots, a_{p(n)}$ with
the following properties:
\roster
\item "(1)" $0 \leq a_i < l$ for all $i$.

\item "(2)" If we set
$s = (a_1 + l k_1) e_1 + \cdots + (a_{p(n)} + l k_{p(n)}) e_{p(n)}$
for any integers $k_1, \ldots, k_{p(n)}$, then $s(x_i) \not= 0$ for all $i$.
\endroster
\endproclaim

\demo{Proof}
Let $Z_i$ be the Zariski closure of $\{ x_i \}$.
We pick up closed a point $y_i$ of $Z_i$ such that
$y_i \not\in Z_j$ for all $j \not= i$.
Let $m_i$ be the maximal ideal at $y_i$.
Since $L$ is $f$-ample, there is a positive integer $n_0$ such that
$$
 H^1(X, L^n \otimes m_1 \otimes \cdots \otimes m_r) = 0
$$
for all $n \geq n_0$. Thus we have
$$
    H^0(X, L^n) \to \bigoplus_{i=1}^{r} L^n/m_iL^n
$$
is surjective. Hence there is $t \in H^0(X, L^n)$ with
$t(y_i) \not= 0$ for all $i$.
Since $l H^0(X, L^n) \subseteq \Ker( H^0(X, L^n) \to L^n/m_iL^n )$
for all $i$, $(t + l k)(y_i) \not= 0$ for all $k \in H^0(X, L^n)$ and all $i$.
Therefore, there are integers $a_1, \ldots, a_{p(n)}$ such that
$0 \leq a_i < l$ for all $i$ and
$a_1 e_1 + \cdots + a_{p(n)} e_{p(n)}$ does not vanish at $y_i$ for all $i$.
Thus, if we set
$s = (a_1 + l k_1)e_1 + \cdots + (a_{p(n)} + l k_{p(n)})e_{p(n)}$,
then $s(y_i) \not= 0$ for all $i$.
In particular, $s(x_i) \not= 0$ for all $i$.
\qed
\enddemo

\proclaim{Lemma 5.2}
Let $k$ be a field of characteristic zero, $V$ a vector space over $k$,
and $f$ a polynomial function over $V$ at most degree $d$, that is,
$f \in \bigoplus_{0 \leq i \leq d} \Sym^i(V^{\vee})$.
Let $e_1, \ldots, e_n$ be generators of $V$,
$a_1, \ldots, a_n \in k$, and $c \in k \setminus \{ 0 \}$.
If
$$
    f((a_1 + c i_1)e_1 + \cdots + (a_n + c i_n)e_n) = 0
$$
for all non-negative integers $i_1, \ldots, i_n$ with
$i_1 + \cdots + i_n \leq d$, then $f = 0$.
\endproclaim

\demo{Proof}
Clearly, we may assume $e_1, \ldots, e_s$ form a basis of $V$.
We set
$$
  a_1e_1 + \cdots + a_se_s + a_{s+1}e_{s+1} + \cdots + a_ne_n =
  a'_1 e_1 + \cdots + a'_s e_s.
$$
Then, we have
$$
  (a_1 + ci_1) e_1 + \cdots + (a_s + ci_s) e_s +
  a_{s+1}e_{s+1} + \cdots + a_ne_n =
  (a'_1 +ci_1) e_1 + \cdots + (a'_s + ci_s) e_s.
$$
Thus,
$$
    f((a'_1 + c i_1)e_1 + \cdots + (a'_s + c i_s)e_s) = 0
$$
for all non-negative integers $i_1, \ldots, i_s$ with
$i_1 + \cdots + i_s \leq d$.
Therefore, we may assume that $\{e_1, \ldots, e_n \}$ is a basis.

Let $\{ X_1, \cdots, X_n \}$ be the dual basis of $\{ e_1, \cdots, e_n \}$.
Then, $f$ is an element of a polynomial ring $k[X_1, \cdots, X_n]$
such that $\deg(f) \leq d$ and
$$
    f(a_1 + c i_1, \ldots, a_n + c i_n) = 0
$$
for all non-negative integers $i_1, \ldots, i_n$ with
$i_1 + \cdots + i_n \leq d$.
Changing variables by $Y_i = c^{-1}(X_i - a_i)$, we may assume that
$a_1 = \cdots = a_n = 0$ and $c = 1$.

We prove this lemma by induction on $n$. We set
$$
f = a_0 X_n^d + a_1(X_1, \ldots, X_{n-1})X_n^{d-1} + \cdots +
    a_{d}(X_1, \ldots, X_{n-1}).
$$
Since $f(0, \ldots, 0, a) = 0$ for
all non-negative integer $a$ with $0 \leq a \leq d$, we have
$a_0 = 0$.
We fix non-negative integers $i_1, \ldots, i_{n-1}$ with
$i_1 + \cdots + i_{n-1} \leq 1$. Then
$$f(i_1, \ldots, i_{n-1}, a) = 0$$
for all non-negative integers $a$ with $0 \leq a \leq d-1$. Thus we have
$a_1(i_1, \ldots, i_{n-1}) = 0$.
Hence, since $\deg(a_1) \leq 1$,
by hypothesis of induction, we get $a_1 = 0$.
Next we fix non-negative integers $i_1, \ldots, i_{n-1}$ with
$i_1 + \cdots + i_{n-1} \leq 2$. Then
$$f(i_1, \ldots, i_{n-1}, a) = 0$$
for all non-negative integers $a$ with $0 \leq a \leq d-2$. Thus we have
$a_2(i_1, \ldots, i_{n-1}) = 0$.
Hence, since $\deg(a_2) \leq 2$, by hypothesis of induction,
we get $a_2 = 0$.
Continuing the same procedures, we have
$a_0 = a_1 = a_2 = \cdots = a_d = 0$. Therefore, $f = 0$.
\qed
\enddemo

\proclaim{Theorem 5.3}
Let $f : X \to \Spec(\ZZ)$ be an arithmetic variety and $(H, k)$ an
arithmetically ample Hermitian line bundle on $X$.
Let $x_1, x_2, \ldots, x_r$ be distinct points of $X$
such that the residue field $\kappa(x_i)$ at $x_i$ has the positive
characteristic
for every $i$. Let $Z_m$ be a decomposable hypersurface of $|H_{\infty}^m|$.
If there is a polynomial $d(t)$ with $\deg Z_m \leq d(m)$ for
all $m$, then, for a sufficiently large integer $m$,
there is a section $\phi \in H^0(X, H^m)$ with following properties:
\roster
\item "(1)" $\phi(x_i) \not= 0$ for all $x_i$.

\item "(2)" The divisor $\zeros(\phi_{\infty})$ in $|H_{\infty}^m|$
            does not belong to $Z_m$.

\item "(3)" $||\phi||_{\sup} < 1$.
\endroster
\endproclaim

\demo{Proof}
Replacing $H$ by a higher multiple of $H$, we may assume that
$$
\Sym^m(H^0(X, H)) \to H^0(X, H^m)
$$
is surjective for all $m \geq 1$.
Moreover, by \cite{Zh, Corollary ~4.8},
we may assume that there is a basis $\{ \phi_1, \cdots, \phi_n \}$
of $H^0(X, H)$
as a $\ZZ$-module such that $||\phi_i||_{\sup} < 1$ for all $i$.
We set $r = \max_{1 \leq i \leq n} \{ ||\phi_i||_{\sup} \}$.
Since $\Sym^m(H^0(X, H)) \to H^0(X, H^m)$ is surjective,
$$
     \{ \phi_1^{e_1}\cdots\phi_n^{e_n} \}_{
     \Sb e_1 \geq 0, \ldots, e_n \geq 0, \\
         e_1 + \cdots + e_n = m \endSb}
$$
forms generators of $H^0(X, L^m)$.
Thus, by Lemma~5.1, if $m$ is sufficiently large,
there are integers $a_{e_1 \cdots e_n}$ such that
$0 \leq a_{e_1 \cdots e_n} < l$ and if we set
$$
  \phi = \sum (a_{e_1 \cdots e_n} + k_{e_1 \cdots e_n}l)
               \phi_1^{e_1}\cdots\phi_n^{e_n}
$$
for integers $k_{e_1 \cdots e_n}$,
then $\phi(x_i) \not= 0$ for all $i$,
where $l = \operatorname{rad}(
\operatorname{char}(\kappa(x_1))\cdots\operatorname{char}(\kappa(x_r)))$.
On the other hand, by Lemma~5.2,
we can find integer $k_{e_1 \cdots e_n}$ such that
$k_{e_1 \cdots e_n} \geq 0$, $\sum k_{e_1 \cdots e_n} \leq d(m)$ and that,
if $\phi$ is the same as before,
$\zeros(\phi_{\infty}) \not\in Z_m$.
Moreover, it is easy to see that
$$
||\phi||_{\sup} \leq (p(m)(l-1) + d(m)l)r^m,
$$
where $p(m) = \rank \Sym^m(H^0(X, H))$.
Therefore, since $p(m)$ and $d(m)$ are polynomials and $r < 1$,
if $m$ is sufficiently large, $||\phi||_{\sup} < 1$.
\qed
\enddemo

\subhead 6. Proof of Main Theorem
\endsubhead
Let $X$, $d$, $r$, $(H, k)$ and $(E, h)$ be the same as in Main Theorem.

For a Hermitian metric $e$ of $E$, we set
$$
\Delta(e) =
\left\{ \achern{2}{E, e} - \frac{r-1}{2r} \achern{1}{E, e}^2 \right\} \cdot
\achern{1}{H, k}^{d-2}.
$$
It is easy to see that, for a positive smooth function $\rho$ on $X$,
we have $\Delta(\rho e) = \Delta(e)$. Thus, we may assume that
$\det(h)$ is an Einstein-Hermitian metric of $\det(E_{\infty})$.

Let $\{ h_t \}_{0 \leq t < \infty}$ be a unique smooth solution
of the evolution equation
$$
   h_t^{-1} \partial_t(h_t) = K(h_t) - c I
\tag{6.1}
$$
of $E_{\infty}$ with the initial condition $h_0 = h$,
where
$$
c = \frac{2 \pi (d -1) (c_1(E_{\infty}) \cdot H_{\infty}^{d-2})}
         {r (H_{\infty}^{d-1}) }
$$
and $K(h_t)$ is the mean curvature of $(E_{\infty}, h_t)$
(cf. \cite{Ko, Chap.VI, \S6, \S7, \S8}).

Since $(\det h_t)^{-1} \partial_t(\det h_t) = \tr(h_t^{-1} \partial_t h_t)$,
by (6.1), we have
$$
(\det h_t)^{-1} \partial_t(\det h_t) = K(\det h_t) - rc.
$$
Therefore, $\det h_t$ also satisfies the evolution equation of
$\det(E_{\infty})$.
On the other hand, $\det h_0 = \det h$ is Einstein-Hermitian.
Thus, $\det h_t$ is Einstein Hermitian for all $0 \leq t < \infty$
(cf. \cite{Ko, Chap.VI, Proposition~9.1}).
Hence, we have a smooth function $\varphi(t)$ on $[0, \infty)$
with $\det h_t = \varphi(t) \det h$.

\proclaim{Lemma 6.2}
$\Delta(h_t)$ is a monotone decreasing function of $t$.
In particular, $\Delta(h) \geq \Delta(h_t)$ for all $t \geq 0$.
\endproclaim

\demo{Proof}
Let $t$, $t'$ be real number with $0 \leq t \leq t'$.
By the same way as in \cite{Mo, Theorem~6.3},
$$
\align
\Delta(h_t) - \Delta(h_{t'})
& = (\Delta(h_t) - \Delta(h)) - (\Delta(h_{t'}) - \Delta(h)) \\
& = (\Delta(h_t) - \Delta(\root{r}\of{\varphi(t)} h)) -
    (\Delta(h_{t'}) - \Delta(\root{r}\of{\varphi(t')}h)) \\
& = DL(h_t, \root{r}\of{\varphi(t)} h) -
    DL(h_{t'}, \root{r}\of{\varphi(t')} h) \\
& = DL(h_t, h) - DL(h_{t'}, h),
    \qquad \hbox{($\because$ \cite{Ko, Chap.VI, Lemma~3.23})}
\endalign
$$
where $DL$ is the Donaldson Lagrangian.
Since $\{ h_t \}$ is a solution of the evolution equation,
$DL(h_t, h)$ is a monotone decreasing function
(cf. \cite{Ko, Chap.VI, Proposition~9.1}).
Therefore, $\Delta(h_t) \geq \Delta(h_{t'})$.
\qed
\enddemo

Let us start the proof of Main Theorem.
In the case of $\dim X = 2$, our theorem is true by \cite{Mo}.
We will prove it by induction on $\dim X$.
Clearly we may assume that $X$ is normal.
Let $X \to \Spec(R) \to \Spec(\ZZ)$ be the Stein factorization
of $f : X \to \Spec(\ZZ)$.
Let $x_1, \ldots, x_l$ be the generic points
of irreducible components of singular fibers of $X \to \Spec(R)$.
Since $$
      \left\{ c_2(E_{\eta}) - \frac{r-1}{2r} c_1(E_{\eta})^2 \right\} \cdot
      H_{\eta}^{d-3} \geq 0
$$
on the generic fiber $X_{\eta}$ of $X \to \Spec(R)$,
$$
      \left\{ c_2(\rest{E}{X_q}) - \frac{r-1}{2r} c_1(\rest{E}{X_q})^2
      \right\} \cdot (\rest{H}{X_q})^{d-3} \geq 0
\tag{6.3}
$$
for all smooth fibers $X_q$ of $X \to \Spec(R)$.
By Corollary~4.3 and Theorem~5.3, for a sufficiently large integer $m$, there
is
a section $\phi \in H^0(X, H^m)$ with the following properties:
\roster
\item "(1)" $\phi(x_i) \not= 0$ for all $x_i$.

\item "(2)" $\zeros(\phi_{\infty})$ is smooth on $X_{\infty}$.

\item "(3)" $\rest{E_{\infty}}{\zeros(\phi_{\infty})}$ is semistable.

\item "(4)" $||\phi||_{\sup} < 1$.
\endroster
Let $\zeros{\phi} = Y + a_1 F_1 + \cdots + a_s F_s$ be the decomposition
of $\zeros{\phi}$ into irreducible divisors such that
$Y$ is horizontal and $F_i$'s are vertical.
By (1), all $F_i$'s are smooth fibers.
Here, replacing $H$ by $H^m$, we may assume that $m = 1$.
Then, we have
$$
\multline
\Delta(h_t)  =
      \left\{ \achern{2}{\rest{(E, h_t)}{Y}} -
              \frac{r-1}{2r} \achern{1}{\rest{(E, h_t)}{Y} }^2 \right\} \cdot
       \achern{1}{\rest{(H, k)}{Y}}^{d-3} \\
  + \sum_{i=1}^{s}
      a_i \left\{ c_2(\rest{E}{F_i}) -
              \frac{r-1}{2r} c_1(\rest{E}{F_i})^2 \right\} \cdot
       (\rest{H}{F_i})^{d-3} \\
 - \int_{X_{\infty}} \log(\sqrt{k(\phi, \phi)})
         \left\{ c_2(E_{\infty}, h_t) - \frac{r-1}{2r} c_1(E_{\infty}, h_t)^2
         \right\} \cdot
         c_1(H_{\infty}, k)^{d-3}.
\endmultline
$$
Therefore, by hypothesis of induction, we get
$$
\Delta(h_t) \geq - \int_{X_{\infty}} \log(\sqrt{k(\phi, \phi)})
         \left\{ c_2(E_{\infty}, h_t) - \frac{r-1}{2r} c_1(E_{\infty}, h_t)^2
         \right\} \cdot
         c_1(H_{\infty}, k)^{d-3}.
\tag{6.4}
$$
On the other hand, by \cite{Ko, Chap.VI, Proposition~9.1 and
Lemma~10.15},
$$
  \lim_{t \to \infty} \max_{X} | K(h_t) - c I| = 0.
$$
Thus, by the same way as in \cite{Ko, Chap.IV, Theorem~4.7 and
Theorem~5.7}, for any positive number $\epsilon$,
if $t$ is sufficiently large, there is a non-negative function $u$ on
$X_{\infty}$ such that
$$
   \left\{ c_2(E_{\infty}, h_t) - \frac{r-1}{2r} c_1(E_{\infty}, h_t)^2
\right\}
   \cdot c_1(H_{\infty}, k)^{d-3} \geq (u - \epsilon) c_1(H_{\infty}, k)^{d-1}.
$$
Thus
$$
\align
\Delta(h_t) & \geq - \int_{X_{\infty}} \log(\sqrt{k(\phi, \phi)})
                                     (u - \epsilon) c_1(H_{\infty}, k)^{d-1} \\
            & \geq \epsilon \int_{X_{\infty}} \log(\sqrt{k(\phi, \phi)})
                                            c_1(H_{\infty}, k)^{d-1}
\endalign
$$
for a sufficiently large $t$, which implies that
$$
\lim_{t \to \infty} \Delta(h_t)
\geq \epsilon \int_{X_{\infty}} \log(\sqrt{k(\phi, \phi)})
c_1(H_{\infty}, k)^{d-1}.
$$
Hence we have ${\displaystyle \lim_{t \to \infty} \Delta(h_t) \geq 0}$.
Therefore, by Lemma~{6.2}, we get $\Delta(h) \geq 0$.

\medskip
Next, we consider the equality condition.
Clearly, we may assume that $\det(h)$ is Einstein-Hermitian.

\proclaim{Lemma 6.5}
If $\Delta(h) = 0$, then $(E_{\infty}, h)$ is Einstein Hermitian.
\endproclaim

\demo{Proof}
Let $e$ be another Hermitian metric of $E_{\infty}$.
We set $e' =  \root{r}\of{\det(h)/\det(e)} e$ and
$a = \root{r}\of{\det(e)/\det(h)}$.
By the same way as above, we have
$\Delta(e) - \Delta(h) = DL(e', h)$. Thus, $DL(e', h) \geq 0$.
On the other hand,
$DL(e, h) = DL(ae', h) = DL(e', h) + DL(ae', e')$.
Here, since $\det(e')$ is Einstein-Hermitian, by an easy calculation,
we have
$$
DL(ae', e') = \frac{r \sqrt{-1}}{2}
\int_{X_{\infty}} \partial(\log a) \overline{\partial}(\log a)
c_1(H_{\infty}, k)^{d-2} \geq 0.
$$
Therefore, we have $DL(e, h) \geq 0$.
This show us that $DL(e, h)$ has the absolute minimal value at $e = h$.
Hence, $h$ is Einstein-Hermitian.
\qed
\enddemo

By Lemma~6.5, if $\dim X = 2$, our assertion is trivial.
So we may assume that $\dim X \geq 3$.
We take a section $\phi \in H^0(X, H^m)$ as before.
Then, by (6.4), we have
$$
\Delta(h) \geq - \int_{X_{\infty}} \log(\sqrt{k(\phi, \phi)})
         \left\{ c_2(E_{\infty}, h) - \frac{r-1}{2r} c_1(E_{\infty}, h)^2
\right\}
         \cdot c_1(H_{\infty}, k)^{d-3}.
$$
Since $h$ is Einstein Hermitian, by the same way as in
\cite{Ko, Chap. IV, Theorem~4.7}, there is a non-negative function $u$
on $X_{\infty}$ such that
$$
\left\{ c_2(E_{\infty}, h) - \frac{r-1}{2r} c_1(E_{\infty}, h)^2 \right\}
\cdot c_1(H_{\infty}, k)^{d-3}
= u c_1(H_{\infty}, k)^{d-1}
$$
and that $u$ is identically zero if and only if $(E_{\infty}, h)$
is projective flat.
Here we assume that $u$ is not identically zero.
Then, since $\log(\sqrt{k(\phi, \phi)}) < 0$ for all points $x \in X_{\infty}$,
$$
\Delta(h) \geq
\int_{X_{\infty}} - \log(\sqrt{k(\phi, \phi)}) u c_1(H_{\infty}, k)^{d-1} > 0
$$
This is a contradiction.
\qed


\enddocument